\documentstyle[amscd,amssymb,verbatim,12pt]{amsart}

\theoremstyle{plain}
\newtheorem{Prop}{Proposition}[section]
\newtheorem{Thm}[Prop]{Theorem}
\newtheorem{Cor}[Prop]{Corollary}

\newtheorem{Lem}[Prop]{Lemma}

\theoremstyle{definition}
\newtheorem{Def}[Prop]{Definition}

\theoremstyle{remark}
\newtheorem{Rem}[Prop]{Remark}

 \def\la{\langle}
\def\ra{\rangle}

\def\int{\mathop{\roman{int}}}

\def\1{^{-1}}

\def\diam{\text{diam}}

\def\St{\text{St}}

\def\asdim{\text{asdim}}

\def\LSS{{\mathcal LSS}}

\def\FF{{\mathcal F}}

\def\BB{{\mathcal B}}
\def\CC{{\mathcal C}}
\def\UU{{\mathcal U}}

\def\dokaz{{\bf Proof. }}
\numberwithin{equation}{section}

\begin{document}
\title[
An alternative definition of coarse structures]%
   {An alternative definition of coarse structures
}

\author{J.~Dydak}
\address{University of Tennessee, Knoxville, TN 37996, USA}
\email{dydak@@math.utk.edu}
\author{C.~S.~Hoffland}
\address{University of Tennessee, Knoxville, TN 37996, USA}
\email{hoffland@@math.utk.edu}

\date{ May 20, 2006}


\begin{abstract}

John Roe \cite{Roe lectures} introduced coarse structures for arbitrary sets $X$
by considering subsets of $X\times X$. That definition, while natural for analysts,
is a bit more difficult to digest for topologists and geometers. In this paper
we introduce large scale structures on $X$ via the notion of uniformly
bounded families and we show their equivalence to coarse structures on $X$. That way all basic concepts of large scale geometry
(asymptotic dimension, slowly oscillating functions, Higson compactification)
have natural definitions and basic results from metric geometry carry over
to coarse geometry.

\end{abstract}

\maketitle

\medskip
\medskip

\section{Introduction}\label{section Introduction}

The second author gave a series of seminar lectures on coarse structures
at University of Tennessee
based on \cite{Roe lectures} and \cite{Hof}. After that it became apparent there is a need
for another approach to coarse structures, an approach more suitable for geometers
and topologists. This paper is an attempt to do just that.

Recall that the {\it star} $\St(B,\UU)$ of a subset $B$ of $X$ with respect
to a family $\UU$ of subsets of $X$ is the union of those elements of $\UU$
that intersect $B$.
Given two families $\BB$ and $\UU$ of subsets of $X$,
$\St(\BB,\UU)$ is the family $\{\St(B,\UU)\}$, $B\in\BB$,
of all stars of elements of $\BB$ with respect to $\UU$.

\begin{Def}\label{LSStructureDef}
A {\it large scale structure} $\LSS_X$ on a set $X$ is a non-empty set of families $\BB$
of subsets of $X$ (called {\it uniformly bounded} 
or {\it uniformly $\LSS_X$-bounded} once $\LSS_X$ is fixed)
satisfying the following conditions:
\begin{enumerate}
\item $\BB_1\in\LSS_X$ implies $\BB_2\in\LSS_X$ if each element of $\BB_2$
consisting of more than one point
is contained in some element of $\BB_1$.
\item $\BB_1,\BB_2\in\LSS_X$ implies $\St(\BB_1,\BB_2)\in\LSS_X$.
\end{enumerate}
\end{Def}

We think of (2) above as a generalization of the triangle inequality.

The {\it trivial extension} $e(\BB)$ of a family $\BB$ is defined as
$\BB\cup \{\{x\}\}_{x\in X}$. Recall that $\BB$ is a {\it refinement} of $\BB'$
if every element of $\BB$ is contained in some element of $\BB'$.
Thus, the meaning of (1) of \ref{LSStructureDef}
is that if $\BB\in\LSS_X$, then all refinements of $e(\BB)$ also belong to $\LSS_X$.

\begin{Prop}\label{PropsOfLSS}
Any large scale structure $\LSS_X$ on $X$ has the following
properties:
\begin{enumerate}
\item $\BB\in\LSS_X$ if each element of $\BB$ consists of at most one point.
\item $\BB_1,\BB_2\in\LSS_X$ implies $\BB_1\cup\BB_2\in\LSS_X$.
\end{enumerate}
\end{Prop}
\dokaz 1). Pick any $\BB_1\in\LSS_X$ and notice $\BB_2:=\BB$
satisfies (1) of \ref{LSStructureDef}.
\par 2). Let $\BB_i':=e(\BB_i)$ for $i=1,2$.
Observe $\BB_i'\in \LSS_X$. Therefore $\BB_3=\St(\BB_1',\BB_2')\in\LSS_X$
and notice any element of $\BB_1\cup\BB_2$ is contained
in an element of $\BB_3$.
\hfill $\blacksquare$

We have two basic examples of large scale structures induced by other structures
on $X$. The first one deals with metric spaces, so let us point out
there is no need to restrict ourselves to metrics assuming
only finite values. To emphasize that, let us call $d\colon X\times X\to R_+\cup\infty$
an {\it $\infty$-metric} if it satisfies all the regular axioms of a metric
(with the understanding that $x+\infty=\infty$).
Notice that $\infty$-metrics have the advantage over regular metrics in the fact
that one can easily define the {\it disjoint union} $\bigoplus\limits_{s\in S}(X_s,d_s)$
of any family of $\infty$-metric spaces $(X_s,d_s)$. Namely, put $d(x,y)=\infty$
if $x$ and $y$ belong to different spaces $X_s$ and $X_t$ (those are assumed
to be disjoint). Conversely, any $\infty$-metric space $(X,d)$ is the disjoint
union of its finite components $(C,d|C)$ (two elements belong to the same
finite component if $d(x,y) < \infty$).

\begin{Prop}\label{MetricLSS}
Any $\infty$-metric space $(X,d)$ has a natural large scale structure
$\LSS(X,d)$ defined as follows:
\par  $\BB\in \LSS(X,d)$ if and only if there is $M
> 0$ such that all elements of $\BB$ are of diameter at most $M$.
\end{Prop}
\dokaz If $\BB_1 \in \LSS(X, d)$ and for each $B_\beta \in \BB_2$
consisting of more than one point there is a $B_\alpha \in \BB_1$
containing $B_\beta$, then $\diam(B_\beta) \leq \diam(B_\alpha) \leq
M$ for each $B_\beta \in \BB_2$, whence $\BB_2 \in \LSS(X, d)$. If
$\BB_1, \BB_2 \in \LSS(X, d)$ then there are $M_1, M_2 > 0$ such
that $\diam(B_\alpha) \leq M_1$ and $\diam(B_\beta) \leq M_2$ for
all $B_\alpha \in \BB_1, B_\beta \in \BB_2$. Thus for any $B_\alpha
\in \BB_1$, $\diam(\St(B_\alpha, \BB_2)) \leq 2M_2 + M_1$, whence
$\St(\BB_1, \BB_2) \in \LSS(X, d)$. It follows that $\LSS(X, d)$ is
a large scale structure. \hfill $\blacksquare$

One can generalize \ref{MetricLSS} as follows:
Given certain families $\FF$ of positive functions from an $\infty$-metric space $X$ to reals
one can define $LSS(X,\FF)$ by declaring $\BB\in \LSS(X,\FF)$
if and only if there is $f\in \FF$ such that $\BB$ refines the family of balls
$\{B(x,f(x))\}_{x\in X}$.

One family of interest is all $f$ such that $\lim\limits_{x\to\infty}\frac{f(x)}{d(x,x_0)}=0$,
where $x_0$ is a fixed point in a metric space $X$
(if $X$ is an $\infty$-metric space, one needs to look at each finite component separately). That leads to the {\it sublinear large scale
structure} on $X$.

\begin{Prop}\label{GroupLSS}
Any group $(X,\cdot)$ has a natural large scale structure
$\LSS_l(X,\cdot)$ defined as follows:
\par $\BB\in \LSS_l(X,\cdot)$ if and only if there
is a finite subset $F$ of $X$ such that $\BB$ refines the shifts $\{x\cdot F\}_{x\in X}$ of $F$.
\end{Prop}
\dokaz Notice that if $\BB\ne\emptyset$ refines $\{x\cdot F\}_{x\in X}$
for some finite subset $F$ of $X$, then $e(\BB)$ also refines $\{x\cdot F\}_{x\in X}$.
\par
Suppose $\BB_i$ refines $\{x\cdot F_i\}_{x\in X}$ for $i=1,2$,
where $F_1$ and $F_2$ are finite subsets of $X$. We may enlarge $F_2$
and assume it is symmetric ($y\in F_2$ implies $y^{-1}\in F_2$).
\par
Let $F$ be the set of all products $x\cdot y\cdot z$, where $x\in F_1$ and $y,z\in F_2$.
Given $B\in \BB_1$ pick $a\in X$ such that $B\subset a\cdot F_1$.
If $B'\in \BB_2$ and $u\in B\cap B'$, choose $y\in X$ so that $B'\subset y\cdot F_2$.
Thus $u=a\cdot f_1=y\cdot f_2$, where $f_1\in F_1$ and $f_2\in F_2$.
Therefore $y=a\cdot f_1\cdot f_2^{-1}$
and $B'\subset a\cdot F$ proving that $\St(B,\BB_2)\subset a\cdot F$.
\hfill $\blacksquare$

\begin{Rem}\label{LeftRightIssue} 
Notice that any group $(X,\cdot)$ has another natural large scale structure
$\LSS_r(X,\cdot)$ defined as follows:
\par $\BB\in \LSS_r(X,\cdot)$ if and only if there
is a finite subset $F$ of $X$ such that $\BB$ refines the shifts $\{F\cdot x\}_{x\in X}$ of $F$.
\par Clearly, the two structures coincide if $X$ is Abelian. However, they may differ even for
finitely presented virtually Abelian groups. 
\par Consider $X=\la a,t\mid t^2=1 \text{ and } tat=a^2\ra$.
Notice every element of $X$ has unique representation as $t^ua^v$, where $u=0,1$.
If $\LSS_l(X,\cdot)=\LSS_r(X,\cdot)$, then for $E=\{1,t\}$ there is a finite subset $F$ of $X$
such that for each $x\in X$ there is $y\in X$ satisfying $x\cdot E\subset F\cdot y$.
Pick $k\ge 1$ such that all elements of $F$ can be represented as $t^ua^v$ so that $u=0,1$ and $|v|\leq k$. Put $x=ta^{6k}$ and choose $y\in X$ satisfying $x\cdot E\subset F\cdot y$.
There is $c=0,1$ and $i$ so that $x=t^ca^iy$ and $|i|\leq k$.
Also, there is $d=0,1$ and $j$ so that $x\cdot t=t^da^jy$ and $|j|\leq k$.
\par Case 1: $c=1$. Now $y=a^{6k-i}$ and $d=0$, so $y=a^{-j}ta^{6k}t=a^{12k-j}$.
That means $6k-i=12k-j$ and $6k=j-i$ contradicting $|i|,|j|\leq k$.
\par Case 2: $c=0$. Now $y=a^{-i}ta^{6k}=ta^{6k-2i}$ and $d=1$, so
$y=a^{-j}ta^{6k}t=ta^{12k-2j}$. Thus $12k-2j=6k-2i$ and $6k=2j-2i$
contradicting $|i|,|j|\leq k$.
\end{Rem}

To create a large scale structure on a set $X$ all one needs is a family
$\LSS_X'$ satisfying conditions resembling finite additivity and (2) of \ref{LSStructureDef}.

\begin{Prop}\label{GeneratingLSS}
If $\LSS_X'$ is a set of families in $X$ such that $\BB_1,\BB_2\in\LSS_X'$
implies existence of $\BB_3\in\LSS_X'$ such that
 $\BB_1\cup\BB_2\cup\St(\BB_1,\BB_2)$ refines $\BB_3$, then
the family $\LSS_X$ of all refinements of trivial extensions of elements of $\LSS_X'$ forms a large scale structure on $X$.
\end{Prop}
\dokaz It suffices to show that, given $\BB_1,\BB_2\in\LSS_X'$, $\St(e(\BB_1),e(\BB_2))$
is a refinement of the trivial extension $e(\BB_3)$ for some $\BB_3\in\LSS_X'$.
Choose $\BB_3\in\LSS_X'$ so that $\St(\BB_1\cup\BB_2,\BB_1\cup\BB_2)$
refines it.

Given $B\in\BB_1$ notice $\St(B,e(\BB_2))=B\cup\St(B,\BB_2)$
is contained in $\St(B,\BB_1\cup\BB_2)$.
Also, $\St(x,e(\BB_2))$ is either a point or there is $B\in\BB_2$
containing $x$ in which case $\St(x,\BB_2)$ is contained in $\St(B,\BB_1\cup\BB_2)$.
\hfill $\blacksquare$
\begin{Rem}\label{RemOnGenLSS}
The family $\LSS_X$ in \ref{GeneratingLSS} is said to be {\it generated} by $\LSS_X'$.
A good example is the {\it discrete large scale structure} on any set $X$
generated by all $\BB$ such that $\bigcup\BB$ is finite.
\end{Rem}

Let us show the analog of Theorem 2.55 (p.34) in \cite{Roe lectures}.
Notice the simplicity of our proof.

\begin{Thm} \label{MetrizationThm}
Given a large scale structure $\LSS_X$ on a set $X$ the following
conditions are equivalent:
\begin{itemize}
\item[a.] There is an $\infty$-metric $d_X$ on $X$ such that
$\LSS_X=\LSS(X,d_X)$.
\item[b.] $\LSS_X$ is generated by a countable set.
\end{itemize}
\end{Thm}
\dokaz a)$\implies$b) is obvious as any $\LSS(X,d_X)$ is generated
by the family of $i$-balls, $i\ge 1$.
\par a)$\implies$b). Pick a sequence $\BB_i\in\LSS_X$ generating $\LSS_X$.
Without loss of generality we may assume $\St(\BB_i,\BB_i)$
refines $\BB_{i+1}$ for all $i\ge 1$.
Define the $\infty$-metric $d_X$ on $X$ by setting $d_X(x,y)$ (if $x\ne y$)
equal the smallest $i$ such that there is $B\in\BB_i$ containing both $x$ and $y$.
If no such $i$ exists, put $d_X(x,y)=\infty$.
\par To show the triangle inequality notice that
$0 < d_X(x,y)\leq d_X(y,z)\leq i$ implies $d_X(x,z)\leq i+1$
as both $x$ and $z$ belong to $\St(y,\BB_i)$ which is contained
in some $B\in\BB_{i+1}$.
\par Clearly $\LSS_X\subset \LSS(X,d_X)$ (each $\BB_i$ refines
the family of $(i+1)$-balls in $(X,d_X)$).
Also, any family of $r$-balls in $(X,d_X)$ refines $\BB_i$ for all $i > r$.
Thus $\LSS_X=\LSS(X,d_X)$.
\hfill $\blacksquare$

\section{Coarse structures and their relation to large scale structures}\label{section coarse}

Recall that a {\it coarse structure} $\CC$ on $X$ is a family of subsets $E$
(called {\it controlled sets})
of $X\times X$ satisfying the following properties:
\begin{enumerate}
\item The diagonal $\Delta=\{(x,x)\}_{x\in X}$ belongs to $\CC$.
\item $E_1\in\CC$ implies $E_2\in\CC$ for every $E_2\subset E_1$.
\item $E\in\CC$ implies $E^{-1}\in\CC$, where $E^{-1}=\{(y,x)\}_{(x,y)\in E}$.
\item $E_1,E_2\in\CC$ implies $E_1\cup E_2\in\CC$.
\item $E, F\in\CC$ implies $E\circ F\in\CC$, where $E\circ F$
consists of $(x,y)$ such that there is $z\in X$ so that $(x,z)\in E$
and $(z,y)\in F$.
\end{enumerate}

\begin{Def}\label{BofEandEofB}
Given a family $\BB$ of subsets of $X$ define $\Delta(\BB)$
as $\bigcup\limits_{B\in\BB}B\times B$.
Given $E\subset X\times X$ define $\BB(E)$ as the family of
all $B\subset X$ such that $B\times B\subset E$.
\end{Def}

\begin{Lem} \label{BofEcircF}
Suppose $\BB_1, \BB_2$ are collections in $X$.
If $\Delta(\BB_i)\subset E_i$ for $i=1,2$,
then $\Delta(\St(\BB_1,\BB_2))\subset (E_2\circ E_1)\circ E_2$.
\end{Lem}
\dokaz Let $(x, y) \in \Delta(\St(\BB_1, \BB_2))$. Then for some $B
\in \BB_1$ there are $B_x, B_y \in \BB_2$, containing $x$ and $y$
respectively, such that there are $z_x \in B \cap B_x$ and $z_y \in
B \cap B_y$. Then
\begin{gather*}
(x, z_x) \in B_x \times B_x \subset \Delta(\BB_2) \subset E_2, \\
(z_y, y) \in B_y \times B_y \subset \Delta(\BB_2) \subset E_1, \\
(z_x, z_y) \in B \times B \subset \Delta(\BB_1) \subset E_1,
\end{gather*}
so there is a $z_x \in X$ such that $(x, z_x) \in E_2$ and $(z_x,
z_y) \in E_1$, whence $(x, z_y) \in E_2 \circ E_1$. But then there
is also a $z_y \in X$ such that $(z_y, y) \in E_2$, whence $(x, y)
\in (E_2 \circ E_1) \circ E_2$ as required. \hfill $\blacksquare$

\begin{Lem} \label{BofEcircF}
Suppose $\BB_1, \BB_2$ are collections in $X$.
If $E_i\subset \Delta(\BB_i)$ for $i=1,2$,
then $E_1\circ E_2\subset  \Delta(\St(\BB_2,\BB_1\cup\BB_2))$.
\end{Lem}
\dokaz Suppose $(x,y)\in E_1\circ E_2$. There is $z$ such that
$(x,z)\in E_1$ and $(z,y)\in E_2$. Therefore
one has $B_1\in \BB_1$ and $B_2\in \BB_2$
so that $x,z\in B_1$ and $z,y\in B_2$.
Put $B_3=\St(B_2,\BB_1\cup\BB_2)$ and notice $B_1\cup B_2\subset B_3$. Thus $x,y\in B_3$.
\hfill $\blacksquare$

\begin{Prop} \label{LSSInducesCS}
Every large scale structure $\LSS_X$ on $X$ induces a coarse structure $\CC$ on $X$
as follows: A subset $E$ of $X\times X$ is declared controlled if and only if
there is $\BB\in\LSS_X$ such that $E\subset \bigcup\limits_{B\in\BB}B\times B$.
\end{Prop}
\dokaz By the remarks after Definition 1.1, all refinements of
$e(\BB)$, for $\BB \in \LSS_X$, themselves belong to $\LSS_X$, meaning
that $\{\{x\}\}_{x \in X}$ is a member of $\LSS_X$. Thus
\[\Delta \subset \bigcup_{B \in \{\{x\}\}}B \times B = \bigcup_{x \in X}\{x\}
\times \{x\}\] so $\Delta \in \CC$. Let $E_1 \in \CC$, so there is a
$\BB \in \LSS_X$ such that $E_1 \subset \Delta(\BB)$. $E_2 \subset
E_1$ then $E_2 \subset \Delta(\BB)$ also, whence $E_2 \in \CC$. It
is clear that if $E \subset \Delta(\BB)$ then $E^{-1} \subset
\Delta(\BB)$, so $E^{-1} \in \CC$. If $E_1, E_2 \in \CC$ then there
are families $\BB_1, \BB_2 \in \LSS_X$ such that $E_1 \subset
\Delta(\BB_1)$ and $E_2 \subset \Delta(\BB_2)$. But
\begin{align*}
E_1 \cup E_2 \subset \Delta(\BB_1) \cup \Delta(\BB_2) & =
\left(\bigcup_{B \in \BB_1}B \times B\right) \cup \left(\bigcup_{B \in \BB_2}B \times B\right) \\
& = \bigcup_{B \in \BB_1 \cup \BB_2}B \times B \\
& = \Delta(\BB_1 \cup \BB_2)
\end{align*}
and since, by Proposition 1.2 $\BB_1 \cup \BB_2 \in \LSS_X$, it
follows that $E_1 \cup E_2 \in \CC$. Finally, let $E_1, E_2 \in \CC$
and $\BB_1, \BB_2 \in \LSS_X$ again be as above. Then $E_1 \circ E_2
\subset \Delta(\St(\BB_2,\BB_1\cup \BB_2))$, which, since both
$\BB_1\cup \BB_2$ and $\BB_2$ are members of $\LSS_X$, is itself a
member of $\LSS_X$, completing the proof. \hfill $\blacksquare$

\begin{Prop} \label{CSInducesLSS}
Every coarse structure $\CC$ on $X$ induces a large scale structure $\LSS_X$
on $X$ as follows:
$\BB$ is declared uniformly bounded if and only if there is a controlled set $E$
such that $ \bigcup\limits_{B\in\BB}B\times B\subset E$.
\end{Prop}
\dokaz Let $B_1 \in \LSS_X$; then there is a controlled set $E \in
\CC$ such that $\Delta(\BB_1) \subset E$. Suppose $\BB_2$ is a
family of subsets of $X$ such that for each $B_\beta \in \BB_2$
consisting of more than one point there is a $B_\alpha \in \BB_1$
containing $B_\beta$. Then $\Delta(\BB_2) \subset \Delta(\BB_1)\cup\Delta
\subset E\cup\Delta\in\CC$, whence $\BB_2 \in \LSS_X$. Now suppose that $\BB_1, \BB_2
\in \LSS_X$, so there are $E_1, E_2 \in \CC$ such that $\Delta(\BB_1)
\subset E_1$ and $\Delta(\BB_2) \subset E_2$. But $(E_2 \circ E_1)
\circ E_1$ is controlled, and $\Delta(\St(\BB_1, \BB_2)) \subset
(E_2 \circ E_1) \circ E_2$, whence $\St(\BB_1, \BB_2) \in \LSS_X$. It
follows that $\LSS_X$ is indeed a large scale structure. \hfill
$\blacksquare$

\section{Higson functions and Higson compactification}\label{section Higson}

In this section we discuss relation of large scale structures on a topological space $X$ to compactifications of $X$. Our approach is quite different from that of
\cite{Roe lectures} (pp.26--31) for coarse structures and seems simpler.

Given a large scale structure $\LSS_X$ on $X$, a subset $B$ of $X$ is
{\it bounded} if $\{B\}\in\LSS_X$.

A bounded continuous function $f\colon X\to R$ is called {\it Higson}
if for every $\BB\in\LSS_X$ and for every $\epsilon > 0$
there is a bounded subset $U$ of $X$ such that
$|f(x)-f(y)| < \epsilon$ for all $x,y\in B\setminus U$, $B\in\BB$.

If $X$ is a topological space, then using Higson maps one can construct
a compact space $h(X,\LSS_X)$ and a natural map $i\colon X\to h(X,\LSS_X)$.
Namely, first we construct $i\colon X\to \prod\limits_{f}[\inf(f),\sup(f)]$
by sending $x$ to $\{f(x)\}_{f}$, and then we declare $h(X,\LSS_X)$ to be
the closure of $i(X)$ in $\prod\limits_{f}[\inf(f),\sup(f)]$.

It is of interest to investigate cases where $h(X,\LSS_X)$ is a compactification
of $X$ (called {\it Higson compactification} of $(X,\LSS_X)$), i.e. $i\colon X\to i(X)$
is a homeomorphism. Here is the simplest sufficient condition for $h(X,\LSS_X)$
to be a compactification.

\begin{Prop}\label{HigsonCExists}
Suppose $X$ is a Tychonoff space.
If $\LSS_X$ is a large scale structure such that the family
of all open and bounded subsets of $X$ forms a basis of $X$,
then $h(X,\LSS_X)$ is a compactification of $X$.
\end{Prop}
\dokaz It suffices to show that the family of Higson maps $f\colon X\to [0,1]$
separates points from closed sets. Indeed, given $x_0\in X\setminus A$,
where $A$ is closed, we find $U$ open and bounded such that
$x_0\in U\subset X\setminus A$. Any map $f\colon X\to [0,1]$ such that
$f(x_0)=1$ and $f|(X\setminus U)\equiv 0$ is a Higson map.
\hfill $\blacksquare$

In case of locally compact Tychonoff spaces $X$ we are interested
in the {\it Higson corona} $\nu(X,\LSS_X):=h(X,\LSS_X)\setminus X$ of $X$.

\begin{Cor}\label{HigsonCoronaExists}
Suppose $X$ is a locally compact Tychonoff space.
If $\LSS_X$ is a large scale structure such that all compact subsets of $X$ are bounded,
then $h(X,\LSS_X)$ is a compactification of $X$.
\end{Cor}
\dokaz Notice all open and relatively compact sets in $X$ are bounded
and form a basis of $X$.
\hfill $\blacksquare$

Given a compactification $c(X)$ of a locally compact Tychonoff space $X$
we are interested in constructing a large scale structure $\LSS(c(X),X)$ on $X$
satisfying the following two conditions:

\begin{itemize}
\item[a.] The bounded subsets of $X$ are precisely relatively compact
subsets of $X$.
\item[b.] The Higson maps of $\LSS(c(X),X)$ include restrictions $f|X$ of all continuous
maps $f:c(X)\to R$.
\end{itemize}

Notice $\St(K,\BB)$ is bounded for every relatively compact $K$
and every $\BB\in\LSS(c(X),X)$. That leads to the following definition.

\begin{Def}\label{DefOfProperFamily}
A family $\BB$ is {\it proper} if $\St(K,\BB)$ is relatively compact
for all relatively compact $K\subset X$. Notice that every $B\in\BB$
is relatively compact in that case (consider $K$ consisting of a point in $B$).
\end{Def}

Recall $E\subset X\times X$ is {\it proper} provided both $E[K]$ and $E^{-1}[K]$
are relatively compact for all relatively compact $K\subset X$ (see Definition 2.1 on p.21
in \cite{Roe lectures}).

\begin{Lem} \label{EOfKForDeltaB}
If $\BB$ is a family of subsets of $X$, then $\Delta(\BB)[K]=\St(K,\BB)$.
\end{Lem}
\dokaz Recall that $E[K]$ is the set of all $x'$ such that there is $x\in K$
satisfying $(x',x)\in E$. If $E=\Delta(\BB)$ that means precisely
there is $B\in\BB$ such that $x',x\in B$ and $x\in K$, i.e. $x'\in\St(K,\BB)$.
\hfill $\blacksquare$

\begin{Cor} \label{ProperFamImpliesProperE}
$\BB$ is proper if and only if $\Delta(\BB)$ is proper.
\end{Cor}

\begin{Prop}\label{StarsOfProperFamilies}
If $\BB_1$ and $\BB_2$ are two proper families,
then $\St(\BB_1,\BB_2)$ is a proper family.
\end{Prop}
\dokaz Notice $\St(K,\St(\BB_1,\BB_2))\subset
\St(\St(K,\BB_1),\BB_2)\cup \St(\St(K,\BB_2),\BB_1)$
for every $K\subset X$. If $K$ is relatively compact,
so is $\St(\St(K,\BB_1),\BB_2)\cup \St(\St(K,\BB_2),\BB_1)$.
\hfill $\blacksquare$

A {\it Higson family} relative to compactification $c(X)$
is a proper family $\BB$ satisfying the following property:
For any map $f:c(X)\to R$ and for any $\epsilon > 0$ there is
a relatively compact set $K$ in $X$ such that $|f(x)-f(y)| < \epsilon$
for all $x,y\in B\setminus K$, $B\in\BB$.

\begin{Prop}\label{StarsOfHigsonFamilies}
If $\BB_1$ and $\BB_2$ are two Higson families,
then $\St(\BB_1,\BB_2)$ is a Higson family.
\end{Prop}
\dokaz Suppose $f:c(X)\to R$ is continuous
and $\epsilon > 0$.
Find a relatively compact set $K$ such that
$|f(x)-f(y)| < \epsilon/4$
for all $x,y\in B\setminus K$, $B\in\BB_1$ or $B\in\BB_2$.
Put $L=\St(\St(K,\BB_1),\BB_2))\cup K$.
Suppose $x,y\in \St(B,\BB_2)\setminus L$ for some $B\in\BB_1$
and $|f(x)-f(y)| >\epsilon$. Clearly, both $x$ and $y$ cannot belong to $B$.
We will discuss the case of $x,y\in X\setminus B$, the other cases are similar.
Thus $x\in B_x\in\BB_2$ and $y\in B_y\in\BB_2$
so that there exist $a\in B\cap B_x$ and $b\in B\cap B_y$.
Notice $a,b\in X\setminus K$ (otherwise $x,y\in L$).
Therefore $|f(a)-f(b)| < \epsilon/4$,
$|f(a)-f(x)| < \epsilon/4$, and $|f(y)-f(b)| < \epsilon/4$
resulting in $|f(x)-f(y)| < 3\cdot\epsilon/4$, a contradiction.
\hfill $\blacksquare$

Define $\LSS(c(X),X)$ as consisting of all
Higson families $\BB$.
It is a large scale structure as the trivial extension of a Higson family
is a Higson family and refinements of Higson families are Higson as well.
Notice every continuous $f\colon c(X)\to R$
restricts to  a Higson map $f|X$ of $\LSS(c(X),X)$.
\par
Using Example 2.34 on p.28 in \cite{Roe lectures} consider
the compactification $c(Z)$ of integers ($Z$ is equipped with the discrete
topology) obtained by identifying two different points $u$ and $v$
in the \v Cech-Stone corona $\beta(Z)\setminus Z$. Notice
$\LSS(c(Z),Z)$
is the discrete large scale structure on $Z$ (generated by families
$\BB$ such that $\bigcup\BB$ is finite) whose Higson functions
are all bounded functions $f:Z\to R$, a set larger than restrictions
$f|Z$ of all continuous maps $f\colon c(Z)\to R$.
Thus, the Higson compactification
of $\LSS(c(Z),Z)$ may be larger than $c(X)$.

\section{Asymptotic dimension}\label{section asdim}

Large scale structures offer a very simple definition of
{\it asymptotic dimension}. Namely,
$\asdim(X,\LSS_X)\leq n$ if $\LSS_X$ is generated by families $\BB$ such that
the {\it multiplicity} of $\BB$ is at most $n+1$ (that means each point $x\in X$
is contained in at most $n+1$ elements of $\BB$).

It is well-known that for metric spaces
the condition $\asdim(X)\leq n$ can be expressed by
one of the following equivalent conditions:

\begin{itemize}
\item[a.] for every uniformly bounded family $\BB$ in $X$ there is
a uniformly bounded family $\BB'$ on $X$ of which $\BB$ is a refinement
such that the multiplicity of $\BB'$ is at most $n+1$.
\item[b.] for every $r > 0$ there is
a decomposition of $X$ as $X_0\cup\ldots\cup X_n$
such that the family of $r$-components of each $X_i$
is uniformly bounded.
\end{itemize}

Our first observation is that one can generalize it
to $\infty$-metric spaces without changing the proof.

\begin{Prop} \label{EquivOfTwoDefinitionsOfAsdimMetric}
Suppose $(X,d)$ is an $\infty$-metric space.
If $n\ge 0$, then the following conditions are equivalent:
\begin{itemize}
\item[a.] for every uniformly bounded family $\BB$ in $X$ there is
a uniformly bounded family $\BB'$ on $X$ of which $\BB$ is a refinement
such that the multiplicity of $\BB'$ is at most $n+1$.
\item[b.] for every $r > 0$ there is
a decomposition of $X$ as $X_0\cup\ldots\cup X_n$
such that the family of $r$-components of each $X_i$
is uniformly bounded.
\end{itemize}
\end{Prop}

Actually, the concept of $\asdim(X,d)\leq n$ for $\infty$-metric spaces has the benefit
that one can express what Bell-Dranishnikov \cite{BellDranish A Hurewicz Type}
call $\asdim(X_s,d_s)\leq n$ {\it uniformly} for all $s\in S$
simply by stating $\asdim(\bigoplus\limits_{s\in S}X_s)\leq n$.

We would like to generalize \ref{EquivOfTwoDefinitionsOfAsdimMetric}
to arbitrary large scale structures.
For that we need the notion of {\it $\BB$-components}.
Those are equivalence classes of the relation $x\sim_\BB y$
meaning that there is a finite chain of points $x_0=x,\ldots,x_k=y$
such that for every $i\ge 0$ (and $i\leq k-1$) there is
$B_i\in\BB$ satisfying $x_i,x_{i+1}\in B_i$.

Our generalization of \ref{EquivOfTwoDefinitionsOfAsdimMetric}
has the advantage that its proof is by reduction to \ref{EquivOfTwoDefinitionsOfAsdimMetric}
which shows that the asymptotic dimension of arbitrary
large scale structures can be reduced to asymptotic dimension
of $\infty$-metric spaces.

\begin{Cor} \label{EquivOfTwoDefinitionsOfAsdim}
Suppose $\LSS_X$ is a large scale structure on a set $X$.
If $n\ge 0$, then the following conditions are equivalent:
\begin{itemize}
\item[a.] for every uniformly bounded family $\BB$ in $X$ there is
a uniformly bounded family $\BB'$ on $X$ of which $\BB$ is a refinement
such that the multiplicity of $\BB'$ is at most $n+1$.
\item[b.] for every uniformly bounded family $\BB$ in $X$ there is
a decomposition of $X$ as $X_0\cup\ldots\cup X_n$
such that the family of $\BB$-components of each $X_i$
is uniformly bounded.
\end{itemize}
\end{Cor}
\dokaz a)$\implies$b).
Given $\BB\in\LSS_X$ construct inductively
a sequence of elements $\BB_i\in\LSS_X$ satisfying the following conditions:
\begin{enumerate}
\item $\BB_1=\BB$,
\item $\St(\BB_i,\BB_\i)$ is a refinement of $\BB_{i+1}$ for each $i\ge 1$,
\item the multiplicity of $\BB_i$ is at most $n+1$ for $i > 1$.
\end{enumerate}

Given two points $x,y\in X$ we define $d(x,y)$
as the smallest integer $i$ such that $x,y\in B\in \BB_i$ for some $i$.
If such integer does not exist, we put $d(x,y)=\infty$.

Notice $\asdim(X,d)\leq n$. Therefore one can decompose
$(X,d)$ as $X_0\cup\ldots\cup X_n$
such that the family of $2$-components of each $X_i$
is uniformly bounded by a fixed integer $M$.
That can be translated into $\BB$-components of each $X_i$
being contained in an element of $\BB_{M+1}$.

\par
b)$\implies$a). Given $\BB_1$ put $\BB_2=\St(e(\BB_1),e(\BB_1))$
and find a decomposition of $X$ as $X_0\cup\ldots\cup X_n$
such that the family of $\BB_2$-components of each $X_i$
is uniformly bounded. Consider $\BB_3$ consisting of stars $\St(C,\BB_1)$,
where $C$ is a $\BB_2$-component of some $X_i$.
Clearly, $\BB_1$ refines $\BB_3$, so it remains to show
that the multiplicity of $\BB_3$ is at most $n+1$.
That follows from the observation that $\St(C,\BB_1)\cap \St(C',\BB_1)=
\emptyset$ for every two different $\BB_2$-components $C$ and $C'$ of the same $X_i$
(otherwise $\St(x,\BB_1)$ would interesect both $C$
and $C'$ for any $x\in \St(C,\BB_1)\cap \St(C',\BB_1)$, a contradiction).
\hfill $\blacksquare$

Our final task is to generalize the Hurewicz Theorem for asymptotic dimension of \cite{BellDranish A Hurewicz Type} and \cite{BDLM}.

\par First let us point out that {\it large scale uniform} functions (or {\it bornologous functions}
in terminology of \cite{Roe lectures}) between metric spaces have a very simple
extension to large scale structures: $f\colon (X,\LSS_X)\to (Y,\LSS_Y)$
is large scale uniform if $f(\BB)\in \LSS_Y$ for all $\BB\in\LSS_X$.

\par Given a function $f\colon (X,\LSS_X)\to (Y,\LSS_Y)$ we need to define
the concept of $\asdim(f)\leq n$. Since that has to do with $f^{-1}(\BB)$
for $\BB\in\LSS_Y$, let us define a natural large scale structure
on the disjoint union $\bigoplus\limits_{s\in S}A_s$ for any family $\{A_s\}_{s\in S}$
of subsets of $X$. Since we want the natural projection
$\bigoplus\limits_{s\in S}A_s\to X$ to be large scale uniform,
the natural choice is to call $\BB$ uniformly bounded in $\bigoplus\limits_{s\in S}A_s$
if and only if there is $\CC\in\LSS_X$ such that $\BB|A_s$ refines $\CC$
for all $s\in S$.

Let us adopt the notation of $\bigoplus\BB$ for the disjoint union of any family
$\BB$. Now, $\asdim(f)\leq n$ means that $\asdim(\bigoplus f^{-1}(\BB))\leq n$
for all $\BB\in\LSS_Y$. 

\begin{Thm} \label{HurThm}
If $f\colon (X,\LSS_X)\to (Y,\LSS_Y)$ is a large scale uniform
function, then
$$\asdim(X,\LSS_X)\leq \asdim(f)+\asdim(Y,\LSS_Y).$$
\end{Thm}
\dokaz Let $\asdim(f)=n$ and $\asdim(Y,\LSS_Y)=k$.

Suppose $\BB_1\in\LSS_X$ is a cover. Let us construct by induction
a sequence of covers $\BB_i\in\LSS_X$ and a sequence of covers $\CC_i\in\LSS_Y$
satisfying the following conditions:
\begin{enumerate}
\item $\St(\BB_i,\BB_i)$ refines $\BB_{i+1}$ for all $i\ge 1$.
\item $f(\BB_i)$ refines $\CC_i$.
\item The multiplicity of $\CC_i$ is at most $k+1$.
\item The cover of $\bigoplus f^{-1}(\CC_i)$ induced by $\BB_i$
refines a cover of multiplicity at most $n+1$ that is
a refinement of the cover of $\bigoplus f^{-1}(\CC_i)$ induced by $\BB_{i+1}$.
\item $\St(\CC_i,\CC_i)$ refines $\CC_{i+1}$ for all $i\ge 1$.
\end{enumerate}

Define the $\infty$-metric $d_X$ on $X$ by setting $d_X(x,y)$
equal the smallest $i$ such that there is $B\in\BB_i$ containing both $x$ and $y$.
If no such $i$ exists, put $d_X(x,y)=\infty$.
Create a $\infty$-metric $d_Y$ on $Y$ the same way using the sequence $\CC_i$.
Notice the following properties of $f\colon (X,d_X)\to (Y,d_Y)$:
\begin{itemize}
\item[a.] $\asdim(Y,d_Y)\leq n$.
\item[b.] $\asdim(f)\leq n$.
\item[c.] $f\colon (X,d_X)\to (Y,d_Y)$ is large scale uniform.
\end{itemize}
Indeed, $\LSS(X,d_X)$ is generated by $\BB_i$'s and $\LSS(Y,d_Y)$
is generated by $\CC_i$'s (see the proof of \ref{MetrizationThm}), so a) and c) follow.
Similarly, b) holds.

Since the proof of Hurewicz Theorem in \cite{BDLM} is valid for $\infty$-metric
spaces, one concludes $\asdim(X,d_X)\leq n+k$.
In particular there is a uniformly bounded family $\UU$ in $(X,d_X)$
such that $\BB_1$ refines $\UU$ and the multiplicity of $\UU$ is at most
$k+n+1$. Notice $\UU$ refines $\BB_M$ for some large $M$.
Thus, $\UU\in\LSS_X$ which completes the proof.
\hfill $\blacksquare$

\section{\v Svarc-Milnor Lemma}

\cite{BDM} gives a simple proof of \v Svarc-Milnor Lemma.
Let us use the approach of this paper to offer an explanation of assumptions
in the \v Svarc-Milnor Lemma.

\par Given a function $f\colon X\to Y$ and given a large scale structure
$\LSS_Y$ on $Y$ let us define {\it the induced large scale structure $f^\ast(\LSS_Y)$}
on $X$ as that generated by $f^{-1}(\BB)$, $\BB\in \LSS_Y$.

\begin{Lem} \label{LeftActionsByIso}
If a group $(G,\cdot)$ acts on the left by isometries on a metric space $(X,d)$,
then $\LSS_l(G,\cdot)\subset f^\ast(\LSS(X,d))$ for any $x_0\in X$
and $f(g):=g\cdot x_0$ for $g\in G$.
\end{Lem}
\dokaz Suppose $F\subset G$ is finite. Put $r=\max\{d(x_0,h\cdot x_0)\mid h\in F\}$.
Given $g\in G$ let $U=B(g\cdot x_0,r)$. It suffices to show
$f(g\cdot F)\subset U$. That is obvious as $d(g\cdot h\cdot x_0,g\cdot x_0)=
d(h\cdot x_0,x_0) < r$ as $h\in F$.
\hfill $\blacksquare$

\begin{Lem} \label{LeftActionsByCocompact}
Suppose a group $(G,\cdot)$ acts on the left by isometries on a metric space $(X,d)$, $x_0\in X$
and $f(g):=g\cdot x_0$ for $g\in G$.
$$\LSS_l(G,\cdot)=f^\ast(\LSS(X,d))$$
 if and only if for any bounded
subset $U$ of $G\cdot x_0$ containing $x_0$
the set $\{g\in G\mid (g\cdot U)\cap U\ne\emptyset\}$
is finite.
\end{Lem}
\dokaz In view of \ref{LeftActionsByIso}, we need to analyze $f^\ast(\LSS(X,d))\subset \LSS_l(G,\cdot)$. It holds if and only if, for any $r > 0$,
there is a finite subset $F_r$ of $G$ such that for any $x\in G\cdot x_0$
there is $g_x\in X$ so that $f^{-1}(B(x,r))\subset g_x\cdot F_r$.
\par Put $U=B(x_0,r)$ and assume $F_r=\{g\in G\mid (g\cdot U)\cap U\ne\emptyset\}$
is finite.
If $x=g_x\cdot x_0$, then $f^{-1}(B(x,r))=\{g\in G\mid d(g\cdot x_0,g_x\cdot x_0) < r\}
=\{g\in G\mid g_x^{-1}g\cdot x_0\in B(x_0,r)\}\subset g_x\cdot F_r$.
\par Assume that, for any $r > 0$,
there is a finite subset $F_r$ of $G$ and $g_0\in X$ so that $f^{-1}(B(x_0,r))\subset g_0\cdot F_r$. Consider $U=B(x_0,r)$ (any bounded subset of $G\cdot x_0$
is contained in such ball). If $h\cdot x_0\in (g\cdot U)\cap U$,
then $h\in f^{-1}(B(x_0,r))$, so $h\in g_0\cdot F_r$.
Also, $g^{-1}\cdot h\in g_0\cdot F_r$
which means the set $\{g\in G\mid (g\cdot U)\cap U\ne\emptyset\}$
is finite.
\hfill $\blacksquare$

\end{document}